\documentclass[12pt]{article}

\usepackage{amsmath}

\numberwithin{equation}{section}
\usepackage{amsfonts}
\usepackage{amsthm}
\usepackage{amssymb}
\usepackage{bbm}
\usepackage{centernot}
\usepackage[x11names]{xcolor}
\usepackage{chemfig}
\usepackage{amsmath}
\usepackage{cite}
\usepackage[nodayofweek]{datetime} 
\usepackage{dsfont}
\usepackage{enumitem}
\usepackage{euscript}
\usepackage{faktor}
\usepackage{float}
\usepackage{geometry}
\usepackage{graphicx}
\usepackage{mathrsfs}
\usepackage{mathtools}
\usepackage{polynom}
\usepackage{stmaryrd}
\usepackage{tikz}
\usepackage{tikz-cd}
\usepackage{wasysym}
\usepackage{xfrac}
\usepackage[x11names]{xcolor}
\usepackage{mathtools}
\usepackage{setspace}
\usepackage[all, cmtip]{xy}
\usepackage{comment}
\usepackage{pgfplots}
\pgfplotsset{compat=1.15}
\usepackage{mathrsfs}
\usetikzlibrary{arrows}
\usepackage{thmtools}

\usepackage{hyperref}

\newlist{legal}{enumerate}{10}
\setlist[legal]{label*=\arabic*.}

\newif\ifproofread

\DeclarePairedDelimiter\abs{\lvert}{\rvert}

\newcommand\norm[1]{\left\lVert#1\right\rVert}

\makeatletter
\newcommand*\bigcdot{\mathpalette\bigcdot@{.5}}
\newcommand*\bigcdot@[2]{\mathbin{\vcenter{\hbox{\scalebox{#2}{$\m@th#1\bullet$}}}}}
\makeatother

\usepackage{mathtools}

\newtheorem*{Lemma*}{Lemma}
\newtheorem*{Corollary*}{corollary}
\newtheorem*{theorem*}{Theorem}
\newtheorem*{definition*}{Definition}

\newtheorem{theorem}{Theorem}[section]
\newtheorem{nota}{Notation}[section]
\newtheorem{definition}[theorem]{Definition}
\newtheorem{construction}[theorem]{Construction}
\newtheorem{remark}[theorem]{Remark}
\newtheorem{lemma}[theorem]{Lemma}

\newtheorem{corollary}[theorem]{Corollary}
\newtheorem{claim}[theorem]{Claim}
\newtheorem{conj}[theorem]{Conjecture} 
\newtheorem{question}[theorem]{Open Question} 

\declaretheoremstyle[
  headfont=\normalfont\bfseries,%\itshape,
  numbered=unless unique,
  bodyfont=\normalfont,
  spaceabove=1em plus 0.75em minus 0.25em,
  spacebelow=1em plus 0.75em minus 0.25em,
]{hartending}

\newcommand{\thistheoremname}{}
\newtheorem*{genericthm}{\thistheoremname}

\newcommand{\NN}{\mathbb{N}}

\newcommand{\QQ}{\mathbb{Q}}

\newcommand{\RR}{\mathbb{R}}

\newcommand{\ZZ}{\mathbb{Z}}

\newcommand{\eps}{\epsilon}
\newcommand{\reg}{\operatorname{Reg}}

\newcommand{\Oo}{\mathcal{O}}

\newcommand{\Mm}{\mathcal{M}}

\newcommand{\cO}{\mathcal{O}}

\newcommand{\stab}{\operatorname{stab}}
\newcommand{\diag}{\operatorname{diag}}

\newcommand{\minvec}{(\on{min-vec}_{\eps})_*}

\newcommand{\on}[1]{\operatorname{#1}}

\newcommand{\disp}[1]{\text{disp}\left(#1\right)}

\newcommand{\inn}[1]{\left<#1\right>}

\newcommand{\cov}[1]{\operatorname{cov}(#1)}

\newcommand{\Sp}{\operatorname{span}}

\newcommand{\pid}{\mathrel{\ooalign{$\lneq$\cr\raise.22ex\hbox{$\lhd$}\cr}}}

\AtBeginDocument{}

\title{Directional $p$-Adic Littlewood Conjecture for Algebraic Vectors}
\date{}
\author{
    Yuval Yifrach
}
\Huge
\begin{document}
\maketitle
\begin{abstract}
%Let $X_n$ be the space of unimodular lattices in $\RR^n$ and let $A$ be the full diagonal group in $\on{SL}_n(\RR)$.
%It is known that compact $A$-orbits originate from modules in totally real degree $n$ number fields.
%Our first result shows that for a natural family of compact orbits $(Ax_k)_k$ all originating from a fixed number field $K$, every weak limit of the Haar measures $m_{Ax_k}$ on those orbits must contain the Haar measure $m_{X_n}$ as an ergodic component.
%This result generalizes certain aspects of the work by Aka and Shapira in \cite{Shapira-Aka} to arbitrary dimensions, as well as elements from Shapira-Zheng in \cite{shapira2021translates}.  

For every vector $\overline \alpha\in \RR^n$ and for every rational approximation $(\overline p,q)\in \RR^n\times\RR$ we can associate the displacement vector $q\alpha-\overline p$.
We focus on algebraic vectors, namely $\overline \alpha=(\alpha_1,\dots,\alpha_n)$ such that $1, \alpha_1, \dots, \alpha_n$ span a rank $n$ number field.
For these vectors, we investigate the size of their displacements as well as the distribution of their directions.
We give a new proof to the result of Bugeaud in \cite{YannPAdic} saying that algebraic vectors $\overline \alpha$ satisfy the $p$-adic Littlewood Conjecture. Namely, we prove that
\begin{equation}
    \liminf_{k \to \infty} \left( k \abs{k}_p \right)^{1/n} \| k (\alpha_1, \dots, \alpha_n) \|_\infty = 0.
\end{equation}  
Our new proof lets us classify all limiting distributions, with a special weighting, of the sequence of directions of the defects in the $\varepsilon$-approximations of $(\alpha_1, \dots, \alpha_n)$. Each such limiting measure is expressed as the pushforward of an algebraic measure on $X_n$ to the sphere.  
\end{abstract}
\section{Introduction}

Let $X_n$ denote the space of unimodular lattices in $\RR^n$, commonly identified with the quotient $\on{SL}_n(\RR)/\on{SL}_n(\ZZ)$. Denote $m_{X_n}$ to be the Haar probability measure on $X_n$ coming from the Haar measure on $\on{SL}_n(\RR)$. Let $A \leq \on{SL}_n(\RR)$ denote the full diagonal subgroup.  
A probability measure $\mu$ on $X_n$ is called \textit{algebraic} if there exists a closed subgroup $H \leq \on{SL}_n(\RR)$ such that $\mu$ is $H$-invariant and supported on an $H$-orbit. 
Among the $A$-invariant probability measures on $X_n$, the family of algebraic measures supported on $A$-orbits is significant and well-studied due to its connections with algebraic number theory. 

Our first goal in this paper is to investigate the structure of possible weak limits of algebraic measures supported on $A$-orbits that are related to each other in a specific way. 

Recall that every algebraic ergodic probability measure supported on an $A$-orbit originates from a full module in a number field, as follows:

\begin{construction}\label{const: compact orbit construction}
Let $K$ be a totally real number field of degree $n$, and let $M \subset K$ be a full module. Denote $\sigma_1,\dots,\sigma_n:K\hookrightarrow \RR$ as an ordering of the natural embeddings of $K$. Let $x_M$ be the normalization of the lattice $(\sigma_1,\dots,\sigma_n)(M) \subset \RR^n$ to have co-volume $1$. Then every $x \in X_n$ with a compact $A$-orbit is of the form $x_M$ for some $K, M$ as above. 

Note that if $\mu$ is an algebraic probability measure supported on an $A$-orbit, then this orbit must be compact.
\end{construction}

In \cite{ToriApp}, together with Solan, we demonstrated that Haar measures on compact $A$-orbits exhibit certain non-rigid properties in the following ways:
\begin{enumerate}
    \item We showed (in \cite[Theorem 1.1]{ToriApp}) that weak limits of Haar measures on compact $A$-orbits need not be ergodic. In fact, they can contain any countable collection of ergodic $A$-invariant measures in their ergodic decomposition.
    \item We proved (in \cite[Theorem 1.5]{ToriApp}) that weak limits of Haar measures on compact $A$-orbits need not be probability measures, and any escape of mass can occur.
\end{enumerate}
These results stand in stark contrast to rigidity phenomena observed for unipotent flows. For instance, in \cite{Mozes_Shah_1995}, it was shown that weak limits of ergodic measures invariant under a one-parameter unipotent flow are always ergodic (and invariant) with respect to a subgroup containing this flow.
\begin{question}
A natural question arises for $n\geq 3$: must any nonzero weak limit of Haar measures on compact $A$-orbits (with discriminant tending to infinity) have the Haar measure on $X_n$ as an ergodic component? 
Our construction in \cite{ToriApp} did not account for the full ergodic decomposition of the weak limits. 
Thus, it is possible that $m_{X_n}$ appears as a component in each of our constructed sequences.
\end{question}

This question is related to a more fundamental conjecture of Margulis:
\begin{conj}[Margulis]
	Every $A$-invariant ergodic probability measure on $X_n$ is algebraic if $n\geq 3$.
\end{conj}

The following theorem is a particular case of \cite[Theorem 1.8]{shapira2021translates} and will be one of our main tools.

\begin{theorem}\label{thm: equidist of prime power sublattices}
Let $p$ be a prime. Let $K$ be a totally real number field of degree $n$, and let $M \leq K$ be a lattice (the $\ZZ$-span of a basis for $K$). Fix an ordering $\sigma_1,\dots,\sigma_n$ of the natural embeddings $K\hookrightarrow \RR$, and denote
\begin{equation}
    x_M=\frac{1}{\cov{\sigma(M)}^{1/n}}\sigma(M)\in X_n,
\end{equation}
where $\sigma=(\sigma_1,\dots,\sigma_n)$.

Write $x_M=g\ZZ^n$ for some $g \in \on{SL}_n(\RR)$.

Define, for any $k \in \ZZ$:
\begin{equation}
    a_k=p^{-k\frac{n+1}{n}}\diag(p^k,\dots,p^k,p^{2k})
\end{equation}
and let $x_k=ga_k\ZZ^n$.
Note that since $x_k$ are all normalized sublattices of $x_M$, they have compact $A$-orbits as well.

Then every weak limit $\mu$ of the Haar measures $m_{Ax_k}$ as $k\to \infty$ is algebraic. When $n$ is prime, $m_{Ax_k}\to m_{X_n}$ as $k\to \infty$.
\end{theorem}

Theorem \ref{thm: equidist of prime power sublattices} addresses one aspect of the following question, informed by the above discussion.

\begin{question}\label{ques: rigidity in a fixed number field}
What kind of rigidity should we expect from $A$-invariant ergodic measures coming from a \textbf{fixed} number field? Could such measures exhibit escape of mass? Must they include $m_{X_n}$ as an ergodic component?
\end{question}

\subsection{$p$-Adic Littlewood conjecture for algebraic vectors}

In this section, we discuss an application of Theorem \ref{thm: equidist of prime power sublattices}.

The $p$-Adic Littlewood Conjecture can be stated as follows:

\begin{conj}\label{conj: padic littlewood}
    Let $\alpha \in \mathbb{R}$ and let $p$ be a prime. Then,
    \begin{equation}\label{eq: padic approximation}
        \liminf_{n \to \infty}  \inn{n\alpha} n \abs{n}_p  = 0,
    \end{equation}
    where $\inn{\cdot}$ denotes the fractional part and $\abs{\cdot}_p$ denotes the $p$-adic norm. Equivalently, for any $\epsilon > 0$, there exist $n_\epsilon \to \infty$ as $\epsilon \to 0$ and an integer $m_\epsilon \in \mathbb{N}$ such that
    \begin{equation}
        \left| \alpha - \frac{m_\epsilon}{n_\epsilon} \right| \leq \frac{\epsilon}{n_\epsilon^2 \abs{n_\epsilon}_p}.
    \end{equation}
\end{conj}

This conjecture can be interpreted as an "approximation compromise" in the context of badly approximable numbers. A number $\alpha \in \mathbb{R}$ is called badly approximable if there exists a constant $c > 0$ such that for all $n \in \mathbb{N}$, we have
\begin{equation}
	\inn{n\alpha} n > c.	
\end{equation}
Conjecture \ref{conj: padic littlewood} suggests that multiplying by the $p$-adic norm of $n$ is enough to cause the left hand side in the above inequality to tend to zero on some sequence $(n_k)_k$ satisfying $n_k\rightarrow\infty$.
\begin{remark}
    Einsiedler and Kleinbock proved in \cite{EinsiedlerKleinbock} that the set of exceptions to Conjecture \ref{conj: padic littlewood} has Hausdorff dimension zero.
\end{remark}

For the case where $\alpha$ is a quadratic irrational (i.e., $\mathbb{Q}(\alpha)$ is a degree 2 extension of $\mathbb{Q}$), this conjecture was formulated by de Mathan and Teuli in \cite[Theorem 2.1]{DeMatan}. In fact, a stronger result was proven in \cite{Shapira-Aka}, which we now explain.

Define the one-parameter group $a(t) = \diag(e^t, e^{-t})$ and the lattices
\begin{equation}
	x_{p^n\alpha} := \begin{pmatrix} 1 & p^n\alpha \\ 0 & 1 \end{pmatrix} \mathbb{Z}^2.
\end{equation}
Equation \eqref{eq: padic approximation} is equivalent to the assertion that the $\{a(t)\}_{t \geq 0}$-orbits of $x_{p^n\alpha}$ become unbounded in the space of lattices $X_2$ as $n \to \infty$.
In \cite[Theorem 2.8]{Shapira-Aka}, the authors proved that the $\{a(t)\}_{t \geq 0}$-orbits of $x_{p^n\alpha}$ are not only unbounded but also equidistributed in $X_2$ as $n \to \infty$.

We propose the following conjecture, which seems to be a natural extension of Conjecture \ref{conj: padic littlewood} but lacks a reference in the literature:

\begin{conj}\label{conj: padic for vectors}
    Let $\overline{\alpha} \in \mathbb{R}^n$ be a vector, and let $p$ be a prime. Denote $\inn{\cdot} : \mathbb{R}^n \to [0,1)^n$ to be the $\ell^{\infty}$ distance to the nearest integer vector. Then, we have
    \begin{equation}\label{eq: vec padic approximation}
        \liminf_{k \to \infty} (k \abs{k}_p)^{1/n} \| \inn{k\overline{\alpha}} \|_\infty = 0.
    \end{equation}
\end{conj}

\begin{definition}\label{def: joint algebraic tuple}
    An $n$-tuple of real numbers $\alpha_1, \dots, \alpha_n$ is called a real joint algebraic $n$-tuple if the set $\{1, \alpha_1, \dots, \alpha_n\}$ spans a totally real number field of degree $n+1$. We will often refer to this simply as an "joint algebraic tuple."
\end{definition}

In light of Conjecture \ref{conj: padic for vectors}, and drawing a parallel to the relation between the $p$-Adic Littlewood Conjecture and the result of Aka-Shapira \cite{Shapira-Aka}, Theorem \ref{thm: littlewood for vectors} below provides a positive answer to Conjecture \ref{conj: padic for vectors} for joint algebraic tuples upon taking $\ell$ from the sequence $(p^m)_m$.
\begin{theorem}\label{thm: littlewood for vectors}
    Let $\alpha_1, \dots, \alpha_n$ be a joint algebraic tuple. 
    Then there exists a constant $C$ depending only on $\alpha_1,\dots,\alpha_n$ such that for every $\ell\in \NN$:
	\begin{equation}\label{eq: vec true padic approximation}
		\liminf_{k \to \infty} k^{1/n} \| \inn{k \ell (\alpha_1, \dots, \alpha_n) } \|_\infty\leq \frac{C}{\ell^{1/n}}
	\end{equation}    
    Moreover, the equation
    \begin{equation}\label{eq: vec padic approximation}
        \liminf_{k \to \infty} (k \abs{k}_p)^{1/n} \| \inn{k (\alpha_1, \dots, \alpha_n)} \|_\infty = 0
    \end{equation}
    immediately follows from \eqref{eq: vec true padic approximation} upon taking $\ell=p^m$, thereby confirming that $(\alpha_1, \dots, \alpha_n)$ satisfies Conjecture \ref{conj: padic for vectors}.
\end{theorem}
\begin{remark}[Comparison to \cite{YannPAdic}]
	In \cite[Theorem 1.3]{YannPAdic}, Bugeaud proves that for every joint algebraic tuple $(\alpha_1,\dots,\alpha_n)$ and for every prime $p$,
	\begin{equation}\label{eq: yann padic}
	    \liminf_{k \to \infty} k^{1/n}\abs{k}_p\log k \| \inn{k (\alpha_1, \dots, \alpha_n)} \|_\infty = 0.
	\end{equation}
	Note that in the formulation by Bugeaud, the term $\abs{k}_{p}$ is not taken to power $1/n$ which gives a weaker statement than the complement case.
	However, Bugeaud's proof in \cite{YannPAdic} actually yields the stronger statement:
	\begin{equation}\label{eq: yann strong padic}
		\liminf_{k \to \infty} (k\abs{k}_p)^{1/n}\log k \| \inn{k (\alpha_1, \dots, \alpha_n)} \|_\infty = 0.
	\end{equation}
	Moreover, Bugeaud proves Theorem \ref{thm: littlewood for vectors}, namely
	\begin{equation*}
		\liminf_{k \to \infty} k^{1/n} \| \inn{k \ell (\alpha_1, \dots, \alpha_n) } \|_\infty\leq \frac{C}{\ell^{1/n}}.
	\end{equation*}
	However, the methods in \cite{YannPAdic} and the methods in this paper are different.
	In \cite{YannPAdic}, Bugeaud proves Equation \eqref{eq: vec padic approximation} by analyzing special units in the number field spanned by $1,\alpha_1,\dots,\alpha_n$, while we use dynamical methods via the Dani correspondence.
	Our analysis enables us to study, as reflected in Theorem \ref{thm: directional padic littlewood} below, distributional properties of the directions of displacements of approximating vectors coming from Equation \eqref{eq: vec true padic approximation} for $\ell=p^m$.
\end{remark}
\subsection{Directional $p$-Adic Littlewood}
The methods used to prove Theorem \ref{thm: littlewood for vectors} can be extended to provide more detailed information about the approximations of joint algebraic tuples.

%Before turning to a precise statement, we will state a more basic version of it to motivate the discussion.
%The following theorem says, roughly, that the approximations we see in Theorem \ref{thm: littlewood for vectors} approach the approximated algebraic vector $(\alpha_1,\dots,\alpha_n)$ from every direction.
%\begin{theorem}\label{thm: thin directional}
%	Let $(\alpha_1,\dots,\alpha_n)$ be an algebraic tuple.
%	Then for every open subset of the sphere $U\subset S^{n-1}$ and for every $\eps>0$ there exist infinitely many $k\in \NN$ such that:
%	\begin{enumerate}
%		\item $(k \abs{k}_p)^{1/n} \| \inn{k (\alpha_1, \dots, \alpha_n)} \|_\infty<\eps$;
%		\item The direction of the vector $\inn{k (\alpha_1, \dots, \alpha_n)}\in \RR^n$ belongs to $U$.
%	\end{enumerate}
%\end{theorem}

We will give a more precise and comprehensive statement below in Theorem \ref{thm: directional padic littlewood}.
First, we will introduce the following definitions.
\begin{definition}\label{def: direction of best approximations}
    Given $v \in \mathbb{R}^n$ and $\overline{r} = (\overline{p}, q) \in \mathbb{Z}^n \times \mathbb{N}$, we define the displacement of $v$ by $\overline{r}$ as
    \[
    \disp{\overline{r}, v} = q^{1/n} (qv - \overline{p}),
    \]
    and the normalized displacement of $v$ by $\overline{r}$ as
    \[
    \theta(v, \overline{r}) = \| qv - \overline{p} \|^{-1} (qv - \overline{p}).
    \]
\end{definition}

\begin{definition}\label{rem: should be a definition}
    For every joint algebraic tuple $\overline{\alpha}$ and for every $\eps>0$, we define the set
    \[
    Q_T(\overline{\alpha},\eps) := \left\{ \overline{r} = (\overline{p}, q) \in \mathbb{Z}^n \times \mathbb{N} \text{ primitive} : q < e^{nT}, e^T \| q \overline{\alpha} - \overline{p} \|_\infty < \epsilon \right\}.
    \]
    This set is finite, and we denote it as $Q_T = (\overline{r}_i)_{i=1}^{N_T}$.
    Moreover, for $\overline{r} = (\overline{p}, q) \in \mathbb{Z}^n \times \mathbb{N}$, we define the weight of $\overline{r}$ up to $T$ as
    \[
    w_{\overline{\alpha}}(\overline{r}, T) = \frac{1}{T} \int_{\{t\in [0,T]: \overline r\in Q_t\}} \frac{1}{|Q_t|} dt.
    \]
    This number represents the proportion of $t$'s for which $\overline r\in Q_t$, normalized by the size of each $Q_t$.
\end{definition}

\begin{remark}
    In Definition \ref{rem: should be a definition}, we highlight a sequence of weights $w(\overline{\alpha}, \overline{r})$ associated with each approximation $\overline{r}$. These weights reflect the quality of each approximation and are natural in the statement of Theorem \ref{thm: directional padic littlewood} below. It is also possible to take the uniform counting measure on $Q_T(\overline{\alpha})$, as done in \cite{shapWeiss}. However, our methods, as well as those of \cite{shapWeiss}, do not yield the explicit formulation appearing in Theorem \ref{thm: directional padic littlewood} in this case.
\end{remark}

\begin{nota}\label{not: epsilon}
    For every $\epsilon > 0$ and $n \in \mathbb{N}$, we denote:
    \begin{enumerate}
    	\item $C_{n, \epsilon} = \left\{ (x_1, \dots, x_n, x_{n+1}) \in \mathbb{R}^{n+1} : 0<\| (x_1, \dots, x_n) \|_\infty < \epsilon, \, |x_{n+1}| \leq 1 \right\}$;
    	\item $X_{n+1,\eps}=\{\Lambda\in X_{n+1}:\exists x=(x_1,\dots,x_{n+1})\in \Lambda \text{ s.t. }x\in C_{n,\eps}, (x_1\dots,x_n)\neq \overline 0\}$;
    	\item $\Mm(S^{n-1})$ is the space of positive Borel measures on the sphere $S^{n-1}$.
    \end{enumerate}
\end{nota}

\begin{definition}\label{def: inter and minvec}
	 Given a Borel measure probability measure $\nu$ on $X_n$ and $\eps>0$ we define the distribution of $\eps$-shortest vector to be the probability measure on $S^{n-2}$ defined by:
			\begin{equation}
				\minvec\nu=\int_{X_{n+1,\eps}} \Theta_{\eps}(\Lambda)d\nu(\Lambda)
			\end{equation}
		where $\Theta_{\eps}:X_{n+1,\eps}\rightarrow \Mm(S^{n-2})$ is defined by 
		\begin{equation}
			\Theta_{\eps}(\Lambda)=\frac{1}{\abs{\Lambda\cap C_{n-1,\eps}}}\sum_{v\in \Lambda\cap C_{n-1,\eps}}\delta_{\pi_{\RR^{n-1}}(v)/\norm{\pi_{\RR^{n-1}}(v)}}
		\end{equation}
		and $\pi_{\RR^{n-1}}$ is the orthogonal projection on the first $n-1$ coordinates. 
\end{definition}

\begin{theorem}\label{thm: directional padic littlewood}
	Let $\overline \alpha$ be a joint algebraic tuple. Let $Q_T(\cdot)$ be defined as in Definition \ref{rem: should be a definition}.
	
	For every $k\in \NN$ there exists $\eps_0(k)>0$ such that $\eps_0(k)\rightarrow 0$ as $k\rightarrow\infty$ and the following holds.
	For $\eps>\eps_0(k),T>0$ define $\mu_{k,T,\eps}$ to be the measure
	\begin{equation}
		\mu_{k,T,\eps}=\sum_{\overline r\in Q_T(p^k\overline \alpha,\eps)}w_{p^k\overline \alpha}(\overline r,T)\delta_{\theta(p^k\overline \alpha,\overline r)}.
	\end{equation}
	Then there exists a matrix $U$, independent on $k$, such that:
	\begin{enumerate}
		\item The sequence $(\mu_{k,T,\eps})_T$ converges as $T\rightarrow \infty$ to a measure $\mu_{k,\eps}$. 
		Moreover, there exists $x(k)\in X_n$ such that $Ax(k)$ is compact and $\mu_{k,\eps}=\minvec Um_{Ax(k)}$.
		\item For every weak limit $\mu_{\eps}$ of $(\mu_{k,\eps})_k$ as $k\rightarrow \infty$
	there exists an algebraic $A$-invariant probability measure $\nu$ on $X_n$ such that
	\begin{equation}\label{eq: pk lim}
		\mu_{\eps}=\minvec(\nu(C_{n-1,\eps})^{-1}U\nu\mid_{C_{n-1,\eps}}).
	\end{equation}
	Moreover, $\nu$ is a weak limit of a subsequence of $(m_{Ax(k)})_k$.
	\end{enumerate}
\end{theorem}
\begin{corollary}\label{cor: full support}
	By \cite{einsiedler2006invariant}, there exists $c\in (0,1]$ such that $\nu\geq cm_{X_n}$
	Therefore, by Equation \eqref{eq: pk lim}:
	\begin{equation}
		\mu_{\eps}\geq c\minvec(m_{X_n})
	\end{equation}
	and in particular, for every $\eps>0$ the support of the limiting measure in \eqref{eq: pk lim} is the whole sphere.
\end{corollary}

%\begin{remark}[A higher rank phenomenon]
%The property that every limiting distribution $\mu_{\varepsilon}$ of the weighted displacements is fully supported on the sphere, as stated in Corollary~\ref{cor: full support}, holds for joint algebraic $n$-tuples for any $n$ and generically, as shown in~\cite[Theorem 1.1(4)]{shapWeiss}. This might suggest that the property applies for all $n$ and every vector $\alpha \in \mathbb{R}^n$. However, a counterexample can be constructed for $n = 2$, where a badly approximable number $\alpha$ has a limiting distribution of displacements that is not fully supported on the sphere. Constructing such an $\alpha$ requires additional background and would significantly lengthen this paper. For this reason, we have opted to present the construction in a separate document.
%It is still possible that this property holds for every $\overline{\alpha}\in \RR^n$ for $n\geq 3$ which would make this a higher rank phenomenon.
%\end{remark}

\begin{remark}[Comparison with \cite{shapWeiss}]
Theorem \ref{thm: directional padic littlewood} shares similarities with \cite[Theorem 1.2]{shapWeiss} while also exhibiting notable differences.

Both Theorem \ref{thm: directional padic littlewood}(a) and \cite[Theorem 1.5]{shapWeiss} concern the displacements of approximations of algebraic vectors. Additionally, both results provide information regarding the distribution of these displacement vectors. However, \cite[Theorem 1.5]{shapWeiss} does not give an explicit description of the distribution of the displacement vectors (denoted as $\nu^{(\mathbb{R}^d)}$ in \cite{shapWeiss}). Instead, the authors derive this distribution from a cross-section measure defined using a limiting process.

To obtain an explicit description of this limiting distribution, we employ a different weighting of the sequence of displacements, as defined in Definition \ref{def: direction of best approximations}. Under this approach, the limiting distribution can be explicitly described, as stated in Theorem \ref{thm: directional padic littlewood}(a), as the pushforward of the Haar measure on the periodic orbit $Ax(k)$. 

A limitation of our result is that it focuses solely on the direction of the displacement, without addressing its magnitude. Furthermore, our analysis does not encompass the distributions of the other two invariants discussed in \cite[Theorem 1.2]{shapWeiss} or the distribution related to the best approximations. While we believe our methods could be extended to prove equidistribution for these additional invariants, we have opted to center our discussion around the $p$-adic Littlewood conjecture for algebraic vectors to maintain clarity.

Lastly, our work examines the possible limits of the measures $\mu_{\epsilon, k}$ as $k \to \infty$, rather than restricting attention to a single measure derived from the displacements of approximations of a specific vector.
\end{remark}
\section*{Acknowledgements}
	The author acknowledges the support of ISF grants number 871/17. This work has received funding from
the European Research Council (ERC) under the European Union’s Horizon 2020 research and innovation
programme (Grant Agreement No. 754475).

\section{Notation and Preliminaries}
In this section we introduce some necessary definitions.

\begin{definition}[$O$-notation]\label{def: O notations}
    For two real functions $f,g$ on a set $A$ we write $f\ll g$  if there exists a constant $C$ independent on the parameters of $f$ and $g$ such that $|f|\leq Cg$ on $A$. 
    The notation $O(g)$ will refer to some implicit function $f$ which satisfies $f\ll g$. 
    The notation $\Theta(g)$ will refer to some implicit function $f$ which satisfies $g \ll f\ll g$. 
    Whenever $r$ is a parameter going to $0$ or $\infty$, the notation $o_r(g)$ will refer to some implicit function $f$ which satisfies $f\ll g\cdot h$, for some implicit function $h \to 0$ as $r$ goes to $0$ or $\infty$ respectively. 
\end{definition}

\begin{definition}\label{def: eps cusp}
	For any $\eps>0$ denote $(X_{n})_{<\eps}$ to be the set of lattices that contain a vector $v$ with $\norm{v}<\eps$. Given a lattice $\Lambda$ and $v\in \Lambda$, we say that $v$ is a shortest vector in $\Lambda$ if $v$ minimizes $\{\norm{u}:0\neq u\in \Lambda\}$. If $v$ is unique up to sign, we say that $\Lambda$ has a unique shortest vector. 
\end{definition}

\subsection{Compact Orbits}
Use $\norm{\cdot}$ to denote the $\ell^{\infty}$ norm on $\RR^n$. 
Given a lattice $\Lambda\subset \RR^n$ we use $\cov{\Lambda}$ to denote the co-volume of $\Lambda$.
Let $X_n$ denote the space of unimodular lattices in $\RR^n$ and let $d_{X_n}(\cdot,\cdot)$ %\yynote{this notation is used also for the measures...}
denote the Riemannian metric on $X_n = \on{SL}_n(\RR)/\on{SL}_n(\ZZ)$ coming from the right invariant Riemannian metric $d_{\on{SL}_n(\RR)}(\cdot, \cdot)$. 
Let $\RR^{n-1}_0=\{v\in \RR^n:\sum_iv_i=0\}$. We abuse notations and define $\exp=\exp\circ \diag:\RR^{n-1}_0\to A$ to be the standard parametrization.
We denote by $m_{X_n}$ probability measure on $X_n = \on{SL}_n(\RR)/\on{SL}_n(\ZZ)$ coming from the Haar measure on $\on{SL}_n(\RR)$.
\begin{definition}[Space of Measures]\label{def: space of measures}
	Let $\mathcal M(X_n)$ denote the space of finite measures on $X_n$ endowed with the topology induced by $\mu_k\rightarrow \mu$ if for any $f\in C_c(X_n)$ it holds that $\mu_k(f)\rightarrow \mu(f)$. 
    % We define a metric on $\mathcal M(X_n)$ which induces this topology by letting, for any $\mu_1,\mu_2\in \mathcal M(X_n)$:
	% \begin{equation}
	% 	d(\mu_1,\mu_2)=\sup_{\eps>0} \varepsilon \sup\left\{\abs{\int f\bd\mu_1-\int f \bd\mu_2}:f\in C_c(X_n)\text{ is $1$-Lipschitz and is supported on }\cK_{\eps}\right\}
	% \end{equation}
	% where $\mathcal K_{\eps}=\{x\in X_n: \lambda_1(x) \ge \eps\}$.
    % \osnote{This is not a metric on the right space. This is a Frechet space but not a banach one. 
    % We should mention that here.}
    % \yynote{Are you sure that the weak * topology on probability measures is metrizable? In any case, we are not really using this metric, except for Conjecture \ref{conj: number fields packets equidistributed}. What do you say we give it up?}
    % \osnote{The weak * toplogy on prob measures is metrizable. The weak * topology on finite measures is also metrizable, but without a homogeneous metric. }
\end{definition}

\begin{definition}\label{def: compact orbit}
    For every degree $n$, totally real number field $K$, denote by ${\rm Lat}_K'$ the set of free $\ZZ$-modules of rank $n$ in $K$. We define an equivalence relation on ${\rm Lat}_K'$ by identifying two lattices $\Lambda_1,\Lambda_2 \subset K$ if $\Lambda_1 = k\Lambda_2$ for some $k\in K^\times$. The quotient space is denoted by ${\rm Lat}_K$, and for every $\Lambda\in {\rm Lat}_K'$, denote by $[\Lambda]\in {\rm Lat}_K$ its equivalence class.
    For every rank $n$, $\ZZ$-module $\Lambda\in {\rm Lat}_K'$ consider the lattice $x_\Lambda := \sigma(\Lambda)/(\cov{\sigma(\Lambda)})^{1/n} \in X_n$, where $\sigma_i: K\hookrightarrow\RR;i=1,\dots,n$ is some ordering of the natural embeddings of $K$ and let $\sigma=(\sigma_1,\dots,\sigma_n): K\rightarrow\RR^n$ denote their concatenation. 
    Denote by $\cO_\Lambda = \{k\in K: k\Lambda\subseteq \Lambda\}$. This is a ring. Denote by $\cO_\Lambda^{\times, >0} = \{u\in \cO_\Lambda^\times:\sigma_i(u)>0: i=1,\dots,n\}.$
    For every $U\subseteq \cO_K^{\times, >0}$ denote $A_U = \{\diag(\sigma_1(u), \sigma_2(u),\dots,\sigma_n(u)): u\in U\}$. 
    Note that these definitions depend implicitly on the ordering of the real embeddings of $K$.
\end{definition}

\begin{definition}\label{def: associated lattice}
	Given a joint algebraic tuple $\overline \alpha=(\alpha_1,\dots,\alpha_n)\in \RR^n$, we denote $$\sigma_i:\on{span}_{\QQ}\{1,\alpha_1,\dots,\alpha_n\}\hookrightarrow \RR$$ to be the field embeddings of the number field associated to $\overline \alpha$ ordered in some way. We denote the normalized (namely to an element in $X_n$) image of $\Sp_{\ZZ}\{1,\alpha_1,\dots,\alpha_n\}$ under $(\sigma_1,\dots ,\sigma_{n+1})$ to be $x_{\overline \alpha}$.
\end{definition}

%\begin{definition}[Correspondence between units and integer matrices]\label{def: unit matrix corr}
%	Let $K$ be some totally real rank $n$ number field and let $\sigma_1,\dots,\sigma_n:K\hookrightarrow\RR$ denote some ordering of the natural embeddings of $K$. Given a full module $M$, namely a $\ZZ$-span of a basis for $K$ over $\QQ$, we denote $\Oo_M$ to be the associated order. Namely, $\Oo_M=\{\alpha\in K:\alpha M\subset M\}$. Note that $\Oo_M\leq \Oo_K$ is a subring. The group of units inside $\Oo_M$ is denoted $\Oo_M^{\times}=\Oo_M\cap K^{\times}$. Note that for any $\eps\in \Oo_M^{\times}$, $\eps M=M$. Choose $\mathcal M=\{m_1,\dots,m_n\}$ to be some ordered $\ZZ$-basis for $M$. Then since $\eps M=M$, there exists a unique $\gamma=\gamma_{\eps,\mathcal M}\in \on{SL}_n(\ZZ)$ such that $\eps m_i=\sum_{j=1}^n\gamma_{ij}m_j$ for every $i=1,\dots,n$. 
%\end{definition}

%\begin{definition}[Discriminant and Regulator]
%	Given a full module $M$ in a number field $K$ we denote $\on{Disc}(M)=\cov{\Oo_M}$ and $\reg(M)=\cov{\log \Oo_M^{\times}}$ where $\log \Oo_M^{\times}$ is a lattice in $\RR^{n-1}$ by Dirichlet's unit Theorem.
%\end{definition}

\subsection{Hecke Neighbors}
In this subsection we give some background on the topic of Hecke neighbors.
\begin{definition}[Definition of the $p$-Hecke Neighbors and the Hecke Operator]\label{def: Hecke}
    For every sequence of integers $0\leq k_1\le k_2\le \dots\le k_n$ consider:
    \[a = a_{p; k_1, k_2, \dots, k_n} = \frac{1}{p^{(k_1+\dots+k_n)/n}}\diag(p^{k_1}, p^{k_2}, \dots, p^{k_n})\in \on{SL}_n(\RR).\]
    For every $x=g\on{SL}_n(\ZZ)\in X_n$ denote $T_a(x) = g\on{SL}_n(\ZZ)a \on{SL}_n(\ZZ)$. This set is finite since $a\on{SL}_n(\ZZ)a^{-1}$ is commensurable to $\on{SL}_n(\ZZ)$.
    The size $\#T_a(x) = \#(\on{SL}_n(\ZZ)a\on{SL}_n(\ZZ)/\on{SL}_n(\ZZ))$ depends only on $k_1,\dots,k_n$ and not on $x$.
    Equivalently, 
    $$T_a(x) = \left\{\frac{1}{\sqrt[n]{\cov{x'}}} x': x'\subseteq x \text{ with }x/x' \cong \ZZ/p^{k_1}\ZZ \oplus\cdots \oplus \ZZ/p^{k_n}\ZZ\right\}.$$
    In addition, given a natural number $m$ we define the set of Hecke-neighbors of index $m$ to be
    \begin{equation}
    	T_m(x)=\left\{\frac{1}{\sqrt[n]{\cov{x'}}} x': x'\subseteq x \text{ with }[x':x]=m\right\}.
    \end{equation}
\end{definition} 

\section{Periodic $A$-orbits vs Unipotent Orbits} \label{sec: per vs uni}
In this section we prove a geometric fact about compact $A$-orbits and $A$-orbits of lattices coming from unipotent matrices. 
To state this fact we need the following notation.
\begin{definition}
	\begin{legal}
		\item \label{def: at} We will denote for $t\in \RR$:
			\begin{equation}
				a(t)=\diag(e^t,e^t,\dots, e^t,e^{-nt})\in \on{SL}_{n+1}(\RR).
			\end{equation}
		\item \label{def: a b unipotent} Denote for $\overline \alpha=(\alpha_1,\dots,\alpha_n)\in \RR^n$: 
			\begin{equation}
				u(\overline \alpha)=\begin{pmatrix}
					1& 0& \cdots & \alpha_1\\
					0& 1& \cdots & \alpha_2\\
					\vdots & \vdots & \cdots & \vdots \\
					0& 0& \cdots & 1
				\end{pmatrix}.
			\end{equation}
		\item \label{def: a b unipotent} Given a joint algebraic tuple $\overline \alpha$ we denote an ordering $\sigma_1,\dots,\sigma_{n+1}$ of the natural embeddings of $K=\on{span}_{\QQ}\{1,\alpha_1,\dots,\alpha_n\}$ such that $\sigma_1(\alpha_i)=\alpha_i$ for every $i=1,\dots,n$.
	We denote:
	\begin{equation}
		B(\overline \alpha)=\begin{pmatrix}
				1& \sigma_1(\alpha_1)& \cdots & \sigma_1(\alpha_n)\\
				1& \sigma_2(\alpha_1)& \cdots & \sigma_2(\alpha_n)\\
				\vdots & \vdots & \cdots & \vdots \\
				1& \sigma_{n+1}(\alpha_1)& \cdots & \sigma_{n+1}(\alpha_n)
						\end{pmatrix}.
	\end{equation}
	so that by Definition \ref{def: compact orbit}, $x_{\overline \alpha}=B(\overline \alpha)\ZZ^{n+1}$.
	\end{legal}
\end{definition}
We will prove, roughly speaking, that for joint algebraic tuples $\overline \alpha\in \RR^n$ and for $k\in \NN$, the one parameter orbit $\{a(t)u(p^k\alpha)\ZZ^{n+1}\}_{t>0}$ remains close to $A(B(\overline \alpha)a(T_k)\ZZ^{n+1})$ where $T_k>0$ is a sequence, uniformly in $k$.
\begin{lemma}\label{lem: unipotent close to compact}
	For any joint algebraic tuple $\overline \alpha$ 
	there exists some $U_0\in \on{SL}_{n+1}(\RR)$ such that:
	\begin{equation}
		\lim_{t\rightarrow \infty}\sup_{T_0<0}d(U_0a(t-T_0)B(\overline \alpha)a(T_0),a(t)[a(-T_0)u(\overline \alpha)a(T_0)])=0.
	\end{equation}
\end{lemma}

\begin{proof}
	We start by noting, similarly to \cite[(5.6)]{shapiraAnnals}, that there exist $Q\in \on{GL}_n(\RR)$ and $q_1,\dots,q_{n+1}\in \RR$ such that
	\[
U=\left(
\begin{array}{@{}c|c@{}}
Q & \begin{array}{@{}c@{}} 0 \\ 0 \\ \vdots \end{array} \\
\cline{1-1}
\multicolumn{1}{@{}c}{
  \begin{matrix} q_1 & q_2 & \cdots \end{matrix}
} & q_{n+1}
\end{array}
\right)
\]
	
%	 $q=\begin{pmatrix}
%			q_{11} & q_{12} & 0\\
%			q_{21} & q_{22} & 0\\
%			q_{31} & q_{32} & q_{33}
%	\end{pmatrix}$
satisfies
	\begin{equation}\label{eq: dani translation}
		UB(\overline \alpha)=u(\overline \alpha).
	\end{equation}
	Since the group $a(t)$ expands only the coordinates $(i,n+1)$ for $i=1,\dots,n$, we deduce that $a(t)Ua(-t)\rightarrow U_0$ as $t\rightarrow \infty$ where 
	\begin{equation}
		U_0=\left(
\begin{array}{@{}c|c@{}}
Q & \begin{array}{@{}c@{}} 0 \\ 0 \\ \vdots \end{array} \\
\cline{1-1}
\multicolumn{1}{@{}c}{
  \begin{matrix} 0 & 0 & \cdots \end{matrix}
} & q_{n+1}
\end{array}
\right).
\end{equation}
	
	Conjugating Eq. \eqref{eq: dani translation} by $a(T_0)$ we deduce:
	\begin{equation}
		a(-T_0)UB(\overline \alpha)a(T_0)=a(-T_0)u(\overline \alpha)a(T_0)
	\end{equation}
	so it follows that for any $t>0$:
	\begin{equation}
		[a(t-T_0)Ua(T_0-t)]a(t-T_0)B(\overline \alpha)a(T_0)=a(t)a(-T_0)UA(\overline \alpha)a(T_0)=a(t-T_0)u(\overline \alpha)a(T_0).
	\end{equation}
	Therefore, for every $t$ large enough:
	\begin{equation}
		d(U_0a(t-T_0)B(\overline \alpha)a(T_0),a(t)a(-T_0)u(\overline \alpha)a(T_0))\leq d(U_0,a(t-T_0)Ua(T_0-t))
	\end{equation}
	 which converges to $0$ uniformly, since $T_0<0$ and $a(t)Ua(-t)\rightarrow U_0$ as $t\rightarrow \infty$.
\end{proof}

The following simple claim follows immediately from the definition of $x_k$ and of $a(t)$.
\begin{claim}\label{cl: padic computations}
	The following equations holds for every $\overline \alpha\in \RR^n$ and $\ell\in \NN$:
	\begin{equation}
		x_k=B(\overline \alpha)a(-\frac{k}{n}\log p)\ZZ^{n+1},
	\end{equation}
	\begin{equation}
		u(\ell\overline \alpha)=a\left(\frac{1}{n}\log \ell\right)u(\overline \alpha)a\left(-\frac{1}{n}\log\ell\right).
	\end{equation}
\end{claim}

\begin{remark}\label{rem: the upgrade}
	Lemma \ref{lem: unipotent close to compact} is similar to \cite[Proposition 7.5]{shapWeiss} with an upgrade special for our use. The difference manifests in the following fact. For  each $T_0=-\frac{k}{n}\log p$, the lattices $B(\overline \alpha)a(T_0)\ZZ^{n+1}$ all have compact $A$-orbits. The claim in \cite[Proposition 7.5]{shapWeiss} implies that for each of them separately,
	\begin{equation}
		\lim_{t\rightarrow \infty}d(U_0a(t-T_0)B(\overline \alpha)a(T_0),a(t)[a(-T_0)u(\overline \alpha)a(T_0)])=0.
	\end{equation}
	Using the special relation between the compact orbits $AB(\overline \alpha)a(T_0)\ZZ^{n+1}$, we note that the proof actually gives uniformity of this convergence in $T_0$.
\end{remark}

As explained before the statement of Lemma \ref{lem: unipotent close to compact}, we proved that $a(t)$-orbits of $u(p^k\alpha)\ZZ^{n+1}$ remain within uniformly bounded distance from the corresponding $A$-orbits of $B(\overline \alpha)a(T_k)\ZZ^{n+1}$ for a certain sequence $T_k\rightarrow \infty$. In the following lemma we prove that in fact for every point $y_k$ in the $A$-orbit of $B(\overline \alpha)a(T_k)\ZZ^{n+1}$ there exist infinitely many $t$'s such that $a(t)u(p^k\alpha)\ZZ^{n+1}$ is within bounded distance from $y_k$.
\begin{lemma}\label{lem: prime case}
	For any compact $A$-orbit $Ax\subset X_{n+1}$ and for any Hecke neighbor $x_k\in T_k(x)$, the $a(t)$-orbit of $x_k$ is equidistributed in $Ax_k$.
\end{lemma}
\begin{proof}
	Write $x=g\ZZ^{n+1}$ and let $K_x$ be the number field associated to $Ax$ under the correspondence in Definition \ref{def: compact orbit}.
	Denote by $\Lambda_x$ the lattice in $V_0=(1,\dots,1)^{\perp}$ (not necessarily unimodular) coming from $\stab_A(x)$ under Definition \ref{def: compact orbit}. We claim that the line spanned by $(1,\dots,1,-n)$ is irrational for $\Lambda_x$. Indeed, if this line was rational, this would imply the existence of a unit $\eps\in K_x^{\times}$, such that its associated matrix $A_{\eps}\in \on{SL}_{n+1}(\ZZ)$ satisfies $a(t)g=gA_{\eps}$ for some $t>0$. Therefore, the characteristic polynomial of $A_{\eps}$, and the minimal polynomial of $\eps$, is split over $\RR$ and has $e^t$ as a root of multiplicity $n$. However, by the theory of Galois extensions $\abs{\on{Gal}(\QQ(\eps)/\QQ)}\mid n+1$ so there is no element of multiplicity $n$ as $n$ does not divide $n+1$.
	We deduce that for any $k$, the $a(t)$-orbit of $x_k$ is dense inside $Ax_k$.
	Therefore the line spanned by $(1,\dots,1,-n)$ is equidistributed in $V_0/\Lambda_x$ and therefore the $a(t)$ orbit is equidistributed in $Ax_k$.
\end{proof}

\section{Proofs of the main theorems}\label{sec: proofs}

The first lemma we need is a standard Dani Correspondence which links diophantine approximations of the vector $\overline \alpha$ and geometric properties of an $a(t)$-orbit. We will use the following notation in the lemma.
\begin{nota}
	For every $\eps>0$ and $n\in \NN$ denote
	\begin{equation}
		C_{n,\eps}=\{(x_1,\dots,x_n,x_{n+1})\in \RR^{n+1}:\norm{(x_1,\dots,x_n)}_{\infty}<\eps,\text{ }\abs{x_{n+1}}\leq 1\}.
	\end{equation}
\end{nota}

\begin{lemma}\label{lem: dani}
	Fix some $\overline \alpha\in \RR^n$, $\eps>0$ and $M\in \NN$. Then there exists $0<t\leq \frac{1}{n}\log M$ such that 
	\begin{equation}\label{eq: cone condition}
		a(t)u(\overline \alpha)\ZZ^{n+1}\cap C_{n,\eps}\neq \emptyset
	\end{equation}
	if and only if there exists $\NN\ni m\leq M$ such that 
	\begin{equation}\label{eq: diophantine}
		m^{1/n}\norm{\inn{m\overline \alpha}}<\eps.
	\end{equation}
	Moreover, in the first implication, if $t\rightarrow \infty$ then also $m\rightarrow \infty$.
\end{lemma}
\begin{proof}
	Fix some $\eps,\overline \alpha$ and $M$ as in the statement of the lemma.
	
	Suppose there exists $0<t\leq\frac{1}{n}\log M$ such that Eq. \eqref{eq: cone condition} holds. Denote $\begin{pmatrix}
	m_1\\
	\vdots \\
	m_n \\
	m
\end{pmatrix}\in \ZZ^{n+1}$ to be a vector realizing the validity of Eq. \eqref{eq: cone condition}, namely such that
	\begin{equation}
			\begin{pmatrix}
		e^t(m\alpha_1+m_1)\\
		\vdots \\
		e^t(m\alpha_n+m_n)\\
		e^{-nt}m
	\end{pmatrix}
	\in C_{n,\eps}.
\end{equation}
	Therefore $\norm{e^t\inn{m\overline \alpha}}<\eps$ and $\abs{e^{-nt}m}<1$ so that $e^{nt}>m$ and so $m^{1/n}\norm{\inn{m\overline \alpha}}<\eps$. Moreover, since $t\leq\frac{1}{n}\log M$, it holds that $e^{nt}\leq M$ and so $\abs{e^{-nt}m}<1$ implies that $m\leq M$.
	
	In the reverse direction, suppose Eq. \eqref{eq: diophantine} holds for $m\leq M$. Let $(m_1,\dots,m_n)$ be the vector realizing the distance to the nearest integer in this equation, namely such that
	\begin{equation}
		m^{1/n}\norm{(m\alpha_1-m_1,\dots,m\alpha_n-m_n)}<\eps.
	\end{equation}
	Then we define the vector $\overline u=\begin{pmatrix}
	-m_1\\
	\vdots \\
	-m_n \\
	m
\end{pmatrix}$ and note that for $t=\frac{1}{n}\log m\leq \frac{1}{n}\log M$, it holds that
\begin{equation}
	a(t)u(\overline \alpha)\overline u\in C_{n,\eps}
\end{equation}
so that $a(t)u(\overline \alpha)\ZZ^{n+1}\cap C_{n,\eps}$ and $t\leq \frac{1}{n}\log M$ as desired.

To prove the 'Moreover' part, note that for every fixed $M>0$, $\sup_{i,m\leq M}\inn{m\alpha_i}$ is bounded below away from zero. Therefore, if $t\rightarrow \infty$, $e^t\sup_{i,m\leq M}\inn{m\alpha_i}\rightarrow \infty$ and so 
	\begin{equation}
			\begin{pmatrix}
		e^t(m\alpha_1+m_1)\\
		\vdots \\
		e^t(m\alpha_n+m_n)\\
		e^{-nt}m
	\end{pmatrix}
	\notin C_{n,\eps}.
\end{equation}
\end{proof}

\begin{claim}\label{cl: density}
	Let $\overline \alpha\in \RR^n$ and $M>m>0$. Suppose $\overline m=(m_1,\dots,m_n)$ is such that
	\begin{equation}
		m^{1/n}\norm{(m\alpha_1+m_1,\dots,m\alpha_n+m_n)}_{\infty}<\eps.
	\end{equation}
%	and $\overline m$ is the unique minimizer of the expression $m^{1/n}\norm{(m\alpha_1+m_1,\dots,m\alpha_n+m_n)}$.
	Let $T>0$ and denote 
	\begin{equation}
		\nu_{T}:=(a(\cdot)u(\overline \alpha)\ZZ^{n+1})_*\frac{1}{T}\lambda\mid_{[0,T]}, 
	\end{equation}
	Then for every $(\overline m,m)\in \ZZ^n\times \NN$, writing $\overline \theta=\frac{m\overline \alpha-\overline m}{\norm{m\overline \alpha-\overline m}}$ 
	 we have:
	\begin{equation}
		\minvec\nu_T(\{\overline \theta\})=\sum_{\overline r\in \theta(\overline \alpha, \cdot)^{-1}(\overline \theta)}w_{\overline \alpha}(\overline r,T)
	\end{equation}
	recall that $\theta(\cdot,\cdot)$ is defined by
 \begin{equation}
 	\theta(v,\overline r)=\norm{qv-\overline p}^{-1}(qv-\overline p).
 \end{equation}
\end{claim}
\begin{proof}
	Denote $x_0=u(\overline \alpha)\ZZ^{n+1}.$
	Let $(\overline m',m')\in \ZZ^n\times\NN$, and write $v=a(t)u(\overline \alpha)(\overline m',m')$. Then the property:
	\begin{equation}
		\overline \theta=\pi_{\RR^{n-1}}(v)/\norm{\pi_{\RR^{n-1}}(v)}, \norm{\pi_{\RR^{n-1}}(v)}_{\infty}<\eps,\abs{v_n}\leq 1
	\end{equation}
	is equivalent, by Lemma \ref{lem: dani}, to $(\overline m',m')\in Q_t\cap \theta(\overline \alpha,\cdot)^{-1}(\overline \theta)$.
	Moreover, by definition of $Q_t$ it holds that:
	\begin{equation}
		\abs{Q_t}=\abs{\{v\in a(t)x_0:\norm{\pi_{\RR^{n-1}}(v)}_{\infty}<\eps,\abs{v_n}\leq 1\}}.
	\end{equation}
	Therefore:
	\begin{equation}
		\minvec\nu_T(\{\overline \theta\})=
	\end{equation}
	\begin{equation*}
		=T^{-1}\int_0^T\frac{1}{\abs{\{v\in a(t)x_0:\norm{\pi_{\RR^{n-1}}(v)}_{\infty}<\eps,\abs{v_n}\leq 1\}}}\sum_{v\in a(t)x_0:\norm{\pi_{\RR^{n-1}}(v)}_{\infty}<\eps,\abs{v_n}\leq 1}\delta_{\pi_{\RR^{n-1}}(v)/\norm{\pi_{\RR^{n-1}}(v)}}(\overline \theta)
	\end{equation*}
	\begin{equation*}
		=T^{-1}\int_0^T\frac{1}{\abs{Q_t}}\sum_{\overline r\in \theta(\overline \alpha, \cdot)^{-1}(\overline \theta)}1=\sum_{\overline r\in \theta(\overline \alpha, \cdot)^{-1}(\overline \theta)}w_{\overline \alpha}(\overline r,T)
	\end{equation*}
	as desired.
\end{proof}

\begin{lemma}\label{lem: continuity of minvec}
	If $\nu_k\rightarrow \nu_0$ are probability measures on $X_n$ such that for every $\eps',c>0$
	\begin{equation}\label{eq: nice sets}
		\nu_0(x: x\cap \{\overline v:\norm{(v_1,\dots,v_{n-1})}_{\infty}=\eps'\}\neq \emptyset\text{ or }\abs{v_n}=c)=0.
	\end{equation}
	Then $\minvec\nu_k\rightarrow \minvec\nu_0$.
\end{lemma}

\begin{proof}
	Let $f:S^{n-2}\rightarrow \RR$ be a bounded continuous function and define $\tilde f:X_n\rightarrow \RR$ to be $\tilde f(x)=\int fd\Theta_{\eps}(x)$. The set of discontinuity points for $\tilde f$ is exactly the set
	\begin{equation}
		D=\{x: x\cap \{\overline v:\abs{v_n}=1, \norm{(v_1,\dots,v_{n-1})}_{\infty}=\eps\}\neq \emptyset\}.
	\end{equation}
	Fix $\delta>0$ small. 
	Define the open set $D_{\delta}$ by:
	\begin{equation}
		D_{\delta}=\{x: x\cap \{\overline v:\abs{\abs{v_n}-1}<\delta, \abs{\norm{(v_1,\dots,v_{n-1})}_{\infty}-\eps}<\delta\}\neq \emptyset\}.		
	\end{equation}
	By Eq. \eqref{eq: nice sets} applied twice and by continuity of measures 
	we deduce the following two properties of $D_{\delta}$:
	\begin{enumerate}
		\item $\nu_0(\partial D_{\delta})=0$
		\item There exists $c(\delta)>0$ such that $c(\delta)\rightarrow 0$ as $\delta\rightarrow 0$ and $\nu_0(D_{\delta})<c(\delta)$.
	\end{enumerate}
Since $\nu_k\rightarrow \nu_0$, we can deduce by (b) that for all $l\geq 0$ large enough $\abs{\nu_l(D_{\delta})-\nu_0(D_{\delta})}<\delta$ so $\nu_l(D_{\delta})<\delta+c(\delta)$ for all $l$ large enough. 
	
	Define $\tilde f_{\delta}$ to be a continuous function such that:
	\begin{enumerate}
		\item $\tilde f_{\delta}=\tilde f$ on $D_{\delta}^c$;
		\item $\tilde f_{\delta}\leq \norm{\tilde f}_{\infty}$.
	\end{enumerate}
	This shows that for every $l$ large enough and for $l=0$:
	\begin{equation}
		\abs{\int \tilde f_{\delta}d\nu_{l}-\int \tilde fd\nu_l}\leq(\delta+c(\delta))\norm{\tilde f}_{\infty}.
	\end{equation}
	Since $\tilde f_{\delta}$ is continuous and bounded, we deduce from $\nu_l\rightarrow \nu_0$ weakly, that for all $l$ large enough
	\begin{equation}
		\abs{\int\tilde fd\nu_{l}-\int \tilde fd\nu_0}\leq (\delta+c(\delta))\norm{\tilde f}_{\infty}
	\end{equation}
	which shows that in fact $\int\tilde fd\nu_{l}\rightarrow \int \tilde fd\nu_0$ as $l\rightarrow \infty$. Since $\minvec\nu(f)=\int \tilde fd\nu$, we get the desired claim.
\end{proof}
We leave the following lemma without proof.
\begin{lemma}\label{lem: zero mass}
	Let $Ax_0$ be a compact $A$-orbit in $X_n$ and let $m_{Ax_0}$ be the uniform measure on $Ax_0$. Then for every $\eps',c>0$ and for $\mu\in \{m_{Ax_0},m_{X_n}\}$:
	\begin{equation}
		\mu(x: x\cap \{\overline v:\norm{(v_1,\dots,v_{n-1})}_{\infty}=\eps'\}\neq \emptyset\text{ or }\abs{v_n}=c)=0.
	\end{equation}
\end{lemma}

%We start by proving Theorem \ref{thm: equidist of prime power sublattices}.
%\begin{proof}
%Let $x_k$ denote the sequence of compact orbits coming from the assumption of Theorem \ref{thm: equidist of prime power sublattices}.
%Note that by Definition \ref{def: Hecke} and by definition of $x_k$, it holds that $x_k\in T_{p^{k}}(x)$.
%By Lemma \ref{lem: disc growth} we deduce that $\on{Disc}(x_k)\leq p^{kn^2}\on{Disc}(x)$.
%%	By the claim above we deduce that $\on{disc}(x_k)\geq p^{3k}$ in the setting of Theorem \ref{thm: equidist of prime power sublattices} and $\on{disc}(x_k)\geq k^3$ in the setting of Theorem \ref{thm: equidist of density one sublattices}.
%By Lemma \ref{lem: shapAka lemma}, $\reg(x_k)\ll p^{k}$.	
%We deduce that $\reg(x_k)\geq \on{disc}(x_k)^{1/n^2}$ so by \cite[Corollary 1.7]{ELMV1}, any weak limit of the measures $m_{Ax_k}$ must contain a Haar component as desired.
%\end{proof}

We start by proving Theorem \ref{thm: directional padic littlewood}.
\begin{proof}
	Fix $\overline \alpha$ as in the Theorem, let $\eps(k)>0$ be the minimal such that 
	\begin{equation}
		\liminf_{q\rightarrow \infty}q\inn{p^kq\alpha}<\eps(k).
	\end{equation}
	By Theorem \ref{thm: littlewood for vectors} we can take $\eps(k)=p^{-k/n}$ which therefore satisfies $\eps(k)\rightarrow 0$ as $k\rightarrow\infty$.
	 Denote any choice of weak limit of the measures $\mu_{k,\eps,T}$ as $T\rightarrow \infty$ to be $\mu_{k,\eps}$. 
	
	Define for every $T>0$
	\begin{equation}
		\nu_{T}:=(a(\cdot)u(p^k\overline \alpha)\ZZ^{n+1})_*\frac{1}{T}\lambda\mid_{[0,T]}.
	\end{equation}
	By definition of $\mu_{k,\eps,t}$ and by Claim \ref{cl: density} we know that for every $\overline \theta\in S^{n-2}$:
	\begin{equation}
		\mu_{k,\eps,T}(\{\overline \theta\})=\sum_{\overline r\in \theta(\overline \alpha, \cdot)^{-1}(\overline \theta)}w_{\overline \alpha}(\overline r,T)=\minvec \nu_T(\{\overline \theta\})
	\end{equation}
	which shows, since both $\mu_{k,\eps,T}$ and $\minvec \nu_T$ are supported on finitely many points in $S^{n-2}$, that 
	\begin{equation}
		\mu_{k,\eps,T}=\minvec \nu_T.
	\end{equation}
	By Claim \ref{cl: padic computations}:
	\begin{equation}
		u(p^k\overline \alpha)=a\left(\frac{k}{n}\log p\right)u(\overline \alpha)a\left(-\frac{k}{n}\log p\right).
	\end{equation}
	By Lemma \ref{lem: unipotent close to compact} there exists some matrix $U$ such that
	\begin{equation}
		\lim_{t\rightarrow \infty}\sup_{T_0<0}d(Ua(t-T_0)B(\overline \alpha)a(T_0),a(t)[a(-T_0)u(\overline \alpha)a(T_0)])=0.
	\end{equation}
	In particular, for $T_0=-\frac{k}{n}\log p$ we can invoke Lemma \ref{lem: prime case} to deduce that 
	\begin{equation}
		\nu_T\rightarrow Um_{Ax_{p^k\overline \alpha}}
	\end{equation}
	as $T\rightarrow \infty$.
	By Lemma \ref{lem: zero mass}, 
	\begin{equation}
	m_{Ax_{p^k\overline \alpha}}(	x: \abs{x\cap \{\overline v:\abs{v_n}=1, \norm{(v_1,\dots,v_{n-1})}_{\infty}=\eps\}}\geq 2)=0.	
	\end{equation}
	Thus, Lemma \ref{lem: continuity of minvec} says that 
	\begin{equation}
		\minvec\nu_T\rightarrow \minvec Um_{Ax_{p^k\overline \alpha}}
	\end{equation}
	as $T\rightarrow \infty$, which shows part (a) of the theorem. 
	
	For part (b), we let $\mu_{\eps}$ be a weak limit of $\mu_{k,\eps}$. By part (a) we know that
	\begin{equation}\label{eq: sequence equality}
		\mu_{k,\eps}=\minvec Um_{Ax(k)}.
	\end{equation}
	where $x(k)=x_{p^k\overline \alpha}$. By Theorem \ref{thm: equidist of prime power sublattices}	we know that every weak limit $\nu$ of $m_{Ax(k)}$ is algebraic. 
	
	By Lemma \ref{lem: zero mass} we know that 
	$$\nu(	x: \abs{x\cap \{\overline v:\abs{v_n}=1, \norm{(v_1,\dots,v_{n-1})}_{\infty}=\eps\}}\geq 2)=0$$
	so the conditions of Lemma \ref{lem: continuity of minvec} hold, and we can deduce from it and Eq. \eqref{eq: sequence equality} that 
	\begin{equation}
		\mu_{k,\eps}\rightarrow \minvec U\nu
	\end{equation}
	as $k\rightarrow \infty$, as desired.
\end{proof}

Finally, we prove Theorem \ref{thm: littlewood for vectors}.
\begin{proof}[Proof of Theorem \ref{thm: littlewood for vectors}]
Fix $\ell\in \NN$.
By Claim \ref{cl: padic computations}, we can denote
	\begin{equation}
		x_{\ell}=B(\overline \alpha)a(-\frac{1}{n}\log \ell)\ZZ^{n+1},
	\end{equation}
	\begin{equation}
		u(\ell\alpha)=a\left(\frac{1}{n}\log \ell\right)u(\overline \alpha)a\left(-\frac{1}{n}\log \ell\right).
	\end{equation}

Lemma \ref{lem: unipotent close to compact} applied for $T_0=-\frac{1}{n}\log \ell$ says that 
	\begin{equation}
		\lim_{t\rightarrow \infty}\sup_{T_0<0}d(U_0a(t-T_0)B(\overline \alpha)a(T_0),a(t)[a(-T_0)u(\overline \alpha)a(T_0)])=0
	\end{equation}
so we can deduce that
\begin{equation}\label{eq: prox}
	d(a(t)u(\ell\overline \alpha)\ZZ^{n+1},U_0a(t)x_{\ell})\rightarrow 0
\end{equation}
as $t\rightarrow \infty$.

Moreover, by Lemma \ref{lem: prime case}, $(t\mapsto a(t)x_k)_*\frac{1}{T}m_{[0,T]}$ becomes equidistributed in $Ax_k$ as $t\rightarrow \infty$ according to the Haar measure on the orbit. 
We will need the following simple claim.
\begin{claim}
	The lattice $x_{\ell}$ contains a vector $v$ such that $\norm{v}_{\infty}\leq \ell^{-1/n}$.
	In particular, we will have $x_{\ell}\in X_{n+1,\ell^{-1/n}}$ (recall Notation \ref{not: epsilon}).
\end{claim}
\begin{proof}
	Recall that 	
	\begin{equation}
		B(\overline \alpha)=\begin{pmatrix}
				1& \sigma_1(\alpha_1)& \cdots & \sigma_1(\alpha_n)\\
				1& \sigma_2(\alpha_1)& \cdots & \sigma_2(\alpha_n)\\
				\vdots & \vdots & \cdots & \vdots \\
				1& \sigma_{n+1}(\alpha_1)& \cdots & \sigma_{n+1}(\alpha_n)
						\end{pmatrix}.
	\end{equation}
	and therefore 
	\begin{equation}
		B(\overline \alpha)a(-\frac{1}{n}\log \ell)=\begin{pmatrix}
				\ell^{-1/n}& \ell^{-1/n}\sigma_1(\alpha_1)& \cdots & \ell^{-1}\sigma_1(\alpha_n)\\
				\ell^{-1/n}& \ell^{-1/n}\sigma_2(\alpha_1)& \cdots & \ell^{-1}\sigma_2(\alpha_n)\\
				\vdots & \vdots & \cdots & \vdots \\
				\ell^{-1/n}& \ell^{-1/n}\sigma_{n+1}(\alpha_1)& \cdots & \ell^{-1}\sigma_{n+1}(\alpha_n)
						\end{pmatrix}.
	\end{equation}
	and so the lattice $x_{\ell}$ contains the vector $\ell^{-1/n}(1,\dots,1)$ which has norm $\ell^{-1/n}$.
	This shows that $x_{\ell}\in X_{n+1,\ell^{-1/n}}$ since the projection of $\ell^{-1/n}(1,\dots,1)$ to the first $n$-coordinates is clearly non-zero.
\end{proof}
Since $\{a(t)x_{\ell}\}_{t>0}$ becomes dense in the periodic orbit $Ax_{\ell}$ and by the claim above, we deduce that there exist unboundedly many positive $t$'s such that $a(t)x_{\ell}\in C_{n,\ell^{-1/n}}$.
This implies that there exist unboundedly many positive $t$'s such that $U_0a(t)x_{\ell}\in C_{n,C\ell^{-1/n}}$ for some $C$ depending only on $\overline \alpha$.
By Equation \eqref{eq: prox} we deduce that also $a(t)u(\ell\overline \alpha)\in C_{n,C\ell^{-1/n}}$ for unboundedly many $t>0$.

By the 'Moreover' part of Lemma \ref{lem: dani} this means that for all $k$ large enough,
\begin{equation}
	\liminf_{k\rightarrow \infty}\norm{k^{1/n}\inn{k\ell\overline \alpha}}<\ell^{-1/n}.
\end{equation}
as desired.

\end{proof}

%We need the following claim which is immediate from the definition of the discriminant.
%\begin{claim}\label{cl: discriminant of hecke neighbor}
%	Given $x\in X_3$ with a compact $A$-orbit and $m>0$, denote $x_m$ to be any index $m$-Hecke neighbor. Then $\on{Disc}(x_m)\geq m^3\on{Disc}(x)$.\yynote{verify}
%\end{claim}

\bibliographystyle{plain}
\bibliography{BibErg}{}
\end{document}